# Type **4** is not Computable

Japheth Wood

December 10, 1996

**Abstract**

We extend a recent result of McKenzie, and show that it is an undecidable problem to determine if 4 appears in the typeset of a finitely generated, locally finite variety.

## 1 Introduction

In this paper, we construct a finite algebra $\mathbf{B}(\mathcal{T})$ for any Turing Machine $\mathcal{T}$ such that the type set of the variety generated by $\mathbf{B}(\mathcal{T})$ contains **4** exactly when the Turing Machine $\mathcal{T}$ halts. Thus we will have shown that the problem of deciding whether or not $\mathbf{4} \in \text{typ}\{V(\mathbf{B}(\mathcal{T}))\}$ is recursively equivalent to the Halting Problem.

## 2 What are Turing Machines?

A Turing machine $\mathcal{T}$ is a theoretical model of a computational device that consists of a two-way infinite tape, a read/write head and a finite list of instructions. The language in which the instructions are written consists of the tape symbols $\{0, 1\}$, the motions symbols $\{L, R\}$, and the internal state symbols $\{\mu_i : 0 \leq i < \omega\}$. A tape is a function $t$ from the integers to the set $\{0, 1\}$ of tape symbols.

For our purposes, we'll start the machine $\mathcal{T}$ in what we call the initial configuration, $\mathcal{Q}_0$. That is, with the tape head in state $\mu_1$, reading the $0^{th}$ square of the tape $t_0$, all whose entries are zero. The machine $\mathcal{T}$ proceeds to compute according to its list of instructions. A typical instruction looks like $\mu_i r s T \mu_m$, with $r$ and $s$ tape symbols, $T$ a motion symbol, $1 \leq i < \omega$, and $0 \leq m < \omega$. This means that if the head is in state $\mu_i$, reading an $r$, then it should replace the $r$ with an $s$, move one square to the right or left, depending on whether $T = R$ or $T = L$, and change the machine state to $\mu_m$.

For any given Turing machine, there is a number $n < \omega$ such that for every $i$, $0 < i \leq n$ and every tape symbol $r$, the machine has exactly one command of the form $\mu_i r s T \mu_m$ with $s$ a tape symbol, $T$ a motion symbol and $0 \leq m \leq n$.

There is no command $\mu_0 r s T \mu_m$ beginning with the state $\mu_0$, and we call $\mu_0$ the halting state. If the Turing machine $\mathcal{T}$ ends up in state $\mu_0$ after a finite number of steps, we say that $\mathcal{T}$ halts. Otherwise, $\mathcal{T}$ doesn't halt.

It will be useful to describe a Turing machine configuration by giving its tape $t$, the integer $n$ that the head is visiting, and the state $\mu_i$ that the machine is in. We shall denote this triple $\langle t, n, \mu_i \rangle$ by $\mathcal{Q}$.

Given two Turing configurations $\mathcal{P}, \mathcal{Q}$ and a subset $N$ of the integers, we shall write $\mathcal{P} \leq_N \mathcal{Q}$ to denote that starting with the machine in configuration $\mathcal{Q}$, the Turing machine $\mathcal{T}$ will arrive in configuration $\mathcal{P}$ in finitely many steps, with the machine head only visiting squares in $N$.

We would also like to point out that for a fixed configuration $\mathcal{P}$, the set $\{\mathcal{Q} : \mathcal{P} \leq_N \mathcal{Q}\}$ can be given a natural ordered tree structure with root $\mathcal{P}$ if $\mathcal{P}$ is in the halting state. There is an edge leading from $\mathcal{Q}$ to $\mathcal{Q}'$ iff $\mathcal{Q}' \leq_N \mathcal{Q}$ in one step.

The main theoretical result that we shall use concerning Turing machines is that the Halting Problem is not computable. That is, given a description of a Turing machine, there is no effective procedure to determine whether or not this Turing machine will actually halt.

## 3 The Algebra $\mathbf{B}(\mathcal{T})$

Let $\mathcal{T}$ be a Turing Machine as above with $k$ internal states $\mu_1, \ldots, \mu_k$ and one halting state $\mu_0$. The algebra $\mathbf{B}(\mathcal{T})$ will consist of a universe of $20k+27$ elements and $12k + 10$ basic operations.

### 3.1 The Universe of $\mathbf{B}(\mathcal{T})$

The universe of $\mathbf{B}(\mathcal{T})$ consists of the disjoint union of the bottom element $0$, the set $P = \{P_0, P_1, P_2\}$, the sequential elements $U = \{1, 2, H\}$, and the machine elements $V \cup \overline{V}$ where for notational purposes we put:

$$V_{ir}^s = \{C_{ir}^s, D_{ir}^s, M_i^r\} \text{ for } 0 \leq i \leq k \text{ and } \{r, s\} \subseteq \{0, 1\}$$

$$\overline{V}_{ir}^s = \{\overline{C}_{ir}^s, \overline{D}_{ir}^s, \overline{M}_i^r\} \text{ for } 0 \leq i \leq k \text{ and } \{r, s\} \subseteq \{0, 1\}.$$

$$V_{ir} = V_{ir}^0 \cup V_{ir}^1 \qquad \overline{V}_{ir} = \overline{V}_{ir}^0 \cup \overline{V}_{ir}^1$$

$$V_i = V_{i0} \cup V_{i1} \qquad \overline{V}_i = \overline{V}_{i0} \cup \overline{V}_{i1}$$

and finally,

$$V = \bigcup_i V_i \text{ and } \overline{V} = \bigcup_i \overline{V}_i$$

## 3.2 The Operations of $\mathbf{B}(\mathcal{T})$

There is a partial order $\leq$ on the universe $B(\mathcal{T})$: $x \leq y$ iff either $x = 0$ or $x = y$ or $x = P_2$ and $y \in \{P_0, P_1\}$.

The operation $x \wedge y$ is defined as the greatest lower bound operation for the order $\leq$. Thus, $\langle B(\mathcal{T}), \wedge \rangle$ is a semilattice reduct of $\mathbf{B}(\mathcal{T})$ of height 2.

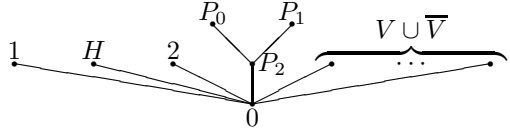

Additionally, we have the three operations $x \wedge^i y$, for $i \in \{0, 1, 2\}$. Put $s_0 = C_{10}^0$, $s_1 = M_1^0$, and $s_2 = D_{10}^0$. Then we take

$$x \wedge^i y = \begin{cases} x & \text{if } x = y \in P \cup V_{10}^0 - \{s_i\} \\ 0 & \text{otherwise.} \end{cases}$$

Now we introduce operations $T_0(x)$ and $T_1(x, y)$ that induce type **4** if $\mathcal{T}$ halts. Note: This is the only change from the type **3** construction in McKenzie [1]

$$T_0(x) = \begin{cases} P_2 & \text{if } x = P_0 \text{ or } P_2 \\ P_0 & \text{if } x = P_1 \\ x & \text{if } x \in V_0 \\ 0 & \text{otherwise.} \end{cases}$$

$$T_1(x, y) = \begin{cases} P_1 & \text{if } x = y = P_0 \text{ or } \{x, y\} = \{P_0, P_1\} \text{ or } \{x, y\} = \{P_0, P_2\} \\ P_2 & \text{if } \{x, y\} \subseteq \{P_1, P_2\} \\ x & \text{if } x = y \in V_{10}^0 \\ 0 & \text{otherwise.} \end{cases}$$

Additionally, we have $J'$, $S_1$ and $S_2$:

$$J'(x, y, z) = \begin{cases} x \wedge z & \text{if } x = y \text{ or } \{x, y\} \subseteq P \\ x & \text{if } x = \overline{y} \in V \cup \overline{V} \\ 0 & \text{otherwise.} \end{cases}$$

$$S_1(u, x, y, z) = \begin{cases} x & \text{if } u \in \{1, 2\} \text{ and } x = y = z \in U \cup V \cup \overline{V} \\ t(x, y, z) & \text{if } \{x, y, z\} \subseteq P \\ 0 & \text{otherwise.} \end{cases}$$

$$S_2(u, v, x, y, z) = \begin{cases} x & \text{if } u = \overline{v} \in U \cup V \cup \overline{V} \text{ and } x = y = z \in V \cup \overline{V} \\ t(x, y, z) & \text{if } \{x, y, z\} \subseteq P \\ 0 & \text{otherwise.} \end{cases}$$

Here, $t(x, y, z)$ is the 3-ary discriminator operation. $t(x, y, z)$ equals $z$ if $x = y$ and $x$ otherwise. Note: except for $J'$, $S_1$ and $S_2$, the element 0 is an **absorbing element** of $B(\mathcal{T})$

## 3.3 The Machine Operations

Finally, we have the machine operations:

For each Turing machine instruction of the form $\mu_i r s L \mu_m$ and for each $t \in \{0, 1\}$ we have the following fundamental operation:

$$L_{irt}(x, y, u) = \begin{cases} u & \text{if } x = y = 1 \text{ and } u \in P \\ C^{s'}_{mt} & \text{if } x = y = 1, u = C^{s'}_{ir} \text{ for some } s' \\ M^t_m & \text{if } x = H, y = 1, u = C^t_{ir} \\ D^s_{mt} & \text{if } x = 2, y = H, u = M^r_i \\ D^{s'}_{mt} & \text{if } x = y = 2, u = D^{s'}_{ir} \text{ for some } s' \\ \overline{v} & \text{if } u \in \overline{V} \text{ and } L_{irt}(x, y, \overline{u}) = v \in V \text{ according to the above.} \\ 0 & \text{otherwise.} \end{cases}$$

$L_{irt}$ simulates the Turing computation at any time that the machine is in state $\mu_i$, reading $r$, and the square to the left of the one being read by the head has $t$ printed on it. We denote this collection of operations as $\mathcal{L}$.

We also have the "reverse" operation:

$$L^o_{irt}(x, y, u) = \begin{cases} v & \text{when } L_{irt}(x, y, v) = u \neq 0. \\ 0 & \text{otherwise.} \end{cases}$$

We denote the collection of these operations as $\mathcal{L}^o$.

For each instruction of the form $\mu_i r s R \mu_m$ in $\mathcal{T}$ and for each $t \in \{0, 1\}$ we have the fundamental operation:

$$R_{irt}(x, y, u) = \begin{cases} u & \text{if } x = y = 1 \text{ and } u \in P \\ C^{s'}_{mt} & \text{if } x = y = 1, u = C^{s'}_{ir} \text{ for some } s' \\ C^s_{mt} & \text{if } x = H, y = 1, u = M^r_i \\ M^t_m & \text{if } x = 2, y = H, u = D^t_{ir} \\ D^{s'}_{mt} & \text{if } x = y = 2, u = D^{s'}_{ir} \text{ for some } s' \\ \overline{v} & \text{if } u \in \overline{V} \text{ and } R_{irt}(x, y, \overline{u}) = v \in V \text{ according to the above.} \\ 0 & \text{otherwise.} \end{cases}$$

These operations are denoted $\mathcal{R}$.

Each of these operations also has its own "reverse" operation.

$$R^o_{irt}(x, y, u) = \begin{cases} v & \text{when } R_{irt}(x, y, v) = u \neq 0. \\ 0 & \text{otherwise.} \end{cases}$$

These operations are denoted by $\mathcal{R}^o$.

There are a total of $4k$ forward and reverse machine operations which we denote by $\mathcal{M}$. Thus, $\mathcal{M} = \mathcal{L} \cup \mathcal{L}^O \cup \mathcal{R} \cup \mathcal{R}^O$. It is also useful to write $\mathcal{M} = \{F_1(x,y,u), \ldots, F_{4k}(x,y,u)\}$, where we have chosen some numbering of these machine operations.

Finally, we have the following set of "error detecting" operations. The purpose of these operations is to detect "bad" configurations of sequential elements. If they are detected, then these operations will produce new elements that will destroy the possibility of a type **4** quotient.

Before we define those, though, we must describe a certain binary relation $\prec$ on the set $U = \{1, H, 2\}$ of sequential elements. This relation is inherent in the way we use these elements. We define $x \prec y$ to be true exactly in the following four cases: $2 \prec 2, 2 \prec H, H \prec 1$, and $1 \prec 1$. Observe that for $1 \leq i \leq 4k$, $F_i(x,y,u)$ is 0 except when $\{x,y\} \subseteq \{1, H, 2\}$ and $x \prec y$.

We now define a set of operations $\mathcal{U}$. For $1 \leq i \leq 4k$, the operations $U_i^1$ and $U_i^2$ are as follows:

$$U_i^1(x,y,z,u) = \begin{cases} \overline{F_i(x,y,u)} & \text{if } x \prec y, x \prec z, y \neq z, F_i(x,y,u) \in V \cup \overline{V}. \\ F_i(x,y,u) & \text{if } x \prec y, y = z. \\ 0 & \text{otherwise.} \end{cases}$$

$$U_i^2(x,y,z,u) = \begin{cases} \overline{F_i(y,z,u)} & \text{if } x \prec z, y \prec z, x \neq y, F_i(y,z,u) \in V \cup \overline{V}. \\ F_i(y,z,u) & \text{if } x = y, y \prec z. \\ 0 & \text{otherwise.} \end{cases}$$

We now define $\mathbf{B}(\mathcal{T})$ to be the algebra with universe $B(\mathcal{T})$ and whose operations are the constant 0, the unary and binary $T_0$ and $T_1$, the operations in $\{\wedge^0, \wedge^1, \wedge^2\} \cup \{\wedge, J', S_1, S_2\}$, the forward and reverse machine operations $\mathcal{M}$, and also all the operations in $\mathcal{U} = \{U_i^1, U_i^2 : 1 \leq i \leq 4k\}$.

## 4   If $\mathcal{T}$ Halts

**Theorem 4.1** *If $\mathcal{T}$ halts, then **4** is in the typeset of $V(\mathbf{B}(\mathcal{T}))$.*

Before proving this, we offer this helpful example:
Example:
Let $\mathcal{T}$ consist of the following commands:

1. $\mu_1 0 1 R \mu_2$
2. $\mu_1 1 0 L \mu_0$
3. $\mu_2 0 1 L \mu_1$
4. $\mu_2 1 0 R \mu_1$

That is, $\mathcal{T}$ has four commands and the three states $\mu_0, \mu_1, \mu_2$. Here is the computation:

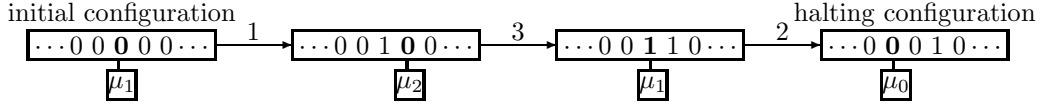

As the machine only actually used three tape squares, we can recover this computation in a third direct power of the algebra $\mathbf{B}(\mathcal{T})$. However, to induce the type **4** quotient, we will use a fourth direct power, the additional coordinate to carry the type **4** lattice structure.

Let us review the algebra $\mathbf{B}(\mathcal{T})$. The universe $B(\mathcal{T})$ consists of the 67 elements $\{0, P_0, P_1, P_2\} \cup \{1, H, 2\} \cup V \cup \overline{V}$. Since $\mathcal{T}$ has three machine states, V consists of the 30 elements:

$$V_0 = \left\{ \begin{array}{ccccc} C_{00}^0 & C_{01}^0 & D_{00}^0 & D_{01}^0 & M_0^0 \\ C_{00}^1 & C_{01}^1 & D_{00}^1 & D_{01}^1 & M_0^1 \end{array} \right\}$$

and

$$(\mathbf{V^0_{10}} \text{ in boldface}) \left\{ \begin{array}{ccccc} \mathbf{C^0_{10}} & C_{11}^0 & \mathbf{D^0_{10}} & D_{11}^0 & \mathbf{M^0_1} \\ C_{10}^1 & C_{11}^1 & D_{10}^1 & D_{11}^1 & M_1^1 \\ C_{20}^0 & C_{21}^0 & D_{20}^0 & D_{21}^0 & M_2^0 \\ C_{20}^1 & C_{21}^1 & D_{20}^1 & D_{21}^1 & M_2^1 \end{array} \right\}$$

$\overline{V}$ is a copy of $V$, with bars on top. Note that we have indicated the elements of $V_0$ and $V_{10}^0$ as they are used to encode the halting and initial configurations, respectively. They appear in the definition of the basic operations $T_0(x)$ and $T_1(x, y)$.

Let $Y = \{-2, -1, 0, 1\}$. These are the indexes of the tape squares visited by the machine head, plus the extra coordinate that we shall use for the type **4** lattice structure. We shall define a subalgebra $\mathbf{K}$ of $\mathbf{A} = \mathbf{B}(\mathcal{T})^\mathbf{Y}$, and show that $\mathbf{K}$ has a type **4** prime quotient.

$a_{-1} = \boxed{\begin{array}{cccc} 1 & H & 2 & 2 \end{array}}$     $a_0 = \boxed{\begin{array}{cccc} 1 & 1 & H & 2 \end{array}}$     $a_1 = \boxed{\begin{array}{cccc} 1 & 1 & 1 & H \end{array}}$

$q^0 = \boxed{\begin{array}{cccc} P_0 & C_{10}^0 & M_1^0 & D_{10}^0 \end{array}}$     $q^1 = \boxed{\begin{array}{cccc} P_1 & C_{10}^0 & M_1^0 & D_{10}^0 \end{array}}$     $q^2 = \boxed{\begin{array}{cccc} P_2 & C_{10}^0 & M_1^0 & D_{10}^0 \end{array}}$

Let $\{a_{-1}, a_0, a_1, q^0, q^1, q^2\}$ be as shown above, and let $\mathbf{K}$ be the subalgebra of $\mathbf{A}$ that they generate.

Let $K_0$ be those elements of $K$ that contain a zero in some coordinate. We will focus on the elements of $K - K_0$, as we will later define a congruence $\theta \in Con\mathbf{K}$ that collapses $K_0$ to a single point.

There are 21 elements of $K - K_0$. Three of them are in $S = \{a_{-1}, a_0, a_1\}$ and the remaining eighteen of them are in $\Lambda = K - K_0 - S$. The elements of $\Lambda$ encode configurations of the Turing computation. The following diagram illustrates six of these – those that arise from the sequential elements in $S$, $q_0$, and the forward and reverse machine operations:

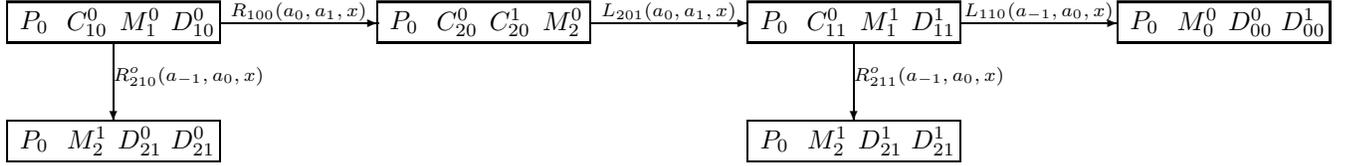

Notice that the Turing machine computation shown on the previous page can be recovered from the top row. The other two elements of $K$ that are shown seem to be unavoidable, as we will require the use of both the forward and reverse machine operations to move freely between elements of $K$ that code initial and halting configurations. However, these extra elements of $K$ do not destroy our result. The astute reader will no doubt be able to calculate the remaining elements of $\Lambda$ with little effort.

We now define $\theta \in ConK$ to be the congruence that collapses $K_0$ to a point. It turns out that an equivalent characterization of $\theta$ is that $(f, g) \in \theta$ iff for all polynomial functions $F(x)$ of $\mathbf{K}$, we have $F(f) = q^2$ iff $F(g) = q^2$. The $\theta$-classes of $\mathbf{K}$ are $K_0$, and the 21 singletons in $\Lambda \cup S$.

We will show later that $\theta$ is strictly meet irreducible in $\mathbf{ConK}$. The unique cover $\overline{\theta}$ of $\theta$ is the congruence generated by $\theta \cup \{(q^0, q^2)\}$.

Now consider the idempotent unary polynomial $e(x) = x \wedge q^0$. The image $e(K)$ consists of $\{q^0, q^2\} \cup T$, where $T$ is some subset of $K_0$. Thus, $e(K)$ consists of two $\overline{\theta}$ classes, one of which splits into two $\theta$ classes. Since $e(K)$ must contain a $\langle \theta, \overline{\theta} \rangle$-minimal set, we must conclude that $N = \{q^0, q^2\}$ is a $\langle \theta, \overline{\theta} \rangle$ trace.

This trace $N$ inherits the semi-lattice meet operation from $K$. It also has a semi-lattice join operation $G^{-1}(T_0(G(T_1(x, y))))|_N$, where $G(x)$ is the composition of the three machine operations (shown in the diagram) that simulate the halting computation of $\mathcal{T}$. $G^{-1}$ is the composition of the corresponding reverse machine operations in reverse order. We also observe that the polynomials of $\mathbf{K}$ that restrict to $N$ are order preserving. Thus we must conclude that $typ\langle \theta, \overline{\theta} \rangle = \mathbf{4}$.

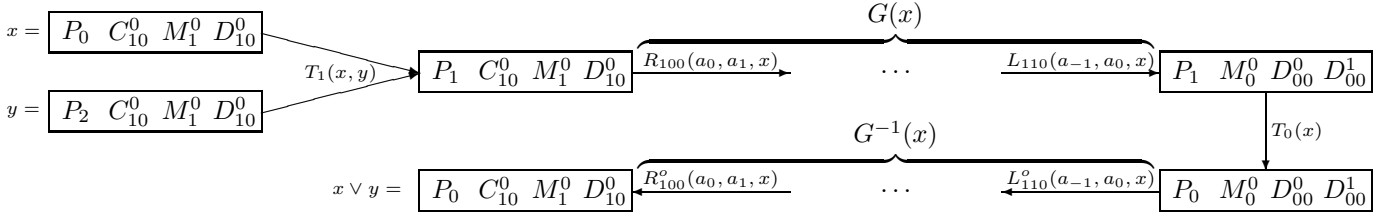

Now we will show the general case.

Assuming that $\mathcal{T}$ halts, let $\mathcal{Q}_0, \mathcal{Q}_1, \ldots, \mathcal{Q}_m$ be the sequence of configurations produced by $\mathcal{T}$, starting with initial configuration $\mathcal{Q}_0 = \langle t_0, 0, \mu_1 \rangle$ where $t_0$ is the zero tape, and for $j \leq m$, $\mathcal{Q}_j = \langle t_j, n_j, \gamma_j \rangle$ and $\gamma_m = \mu_0$.

Let $N = \{n_0, \ldots, n_m\} \cup \{-1, 0, 1\}$. This is essentially the finite interval of integers that the Turing head visits. It is necessary to ensure that $\{-1, 0, 1\} \subseteq N$, as we shall see later that the operations $\wedge^i$ require this. We shall also denote $N$ by $[b, c]$, where $b$ and $c$ are the extreme points of the interval. As in the example above, we shall use some point $s_0 \notin N$ to index the coordinate that houses the type **4** lattice structure. Let $Y = \{s_0\} \cup N$. The following vectors are in $\mathbf{A} = \mathbf{B}(\mathcal{T})^Y$.

For $n \in \{n_0, \ldots, n_m\}$ let $a_n$ be the element of $\mathbf{A}$ such that

$$a_n(x) = \begin{cases} 1 & \text{for } x = s_0 \\ 1 & \text{for } x \in [b, n-1] \\ H & \text{for } x = n \\ 2 & \text{for } x \in [n+1, c] \end{cases}$$

Let $S = \{a_b, \ldots, a_c\}$.

For $i \in \{0, 1, 2\}$ let $q^i \in A$ be the function

$$q^i(x) = \begin{cases} P_i & \text{for } x = s_0 \\ C_{10}^0 & \text{for } x \in [b, -1] \\ M_1^0 & \text{for } x = 0 \\ D_{10}^0 & \text{for } x \in [1, c] \end{cases}$$

**Definition 4.1** *We will work in* $\mathbf{K}$, *the subalgebra of* $\mathbf{A} = \mathbf{B}(\mathcal{T})^Y$ *generated by* $S \cup \{q^0, q^1, q^2\}$.

**Definition 4.2** $\theta$ *is the congruence relation on* $\mathbf{K}$ *such that* $(f, g) \in \theta$ *iff* $\{f, g\} \subseteq K$ *and for all polynomial functions* $F(x)$ *of* $\mathbf{K}$ *we have* $F(f) = q^2$ *iff* $F(g) = q^2$.

**Definition 4.3** $\overline{\theta}$ *is the congruence on* $\mathbf{K}$ *generated by* $\theta \cup \{(q^0, q^2)\}$.

**Definition 4.4** $Cfg_\star$ is the set of all configurations $\mathcal{Q}$ such that $\mathcal{Q}_m \leq_N \mathcal{Q}$.

Recall that $Cfg_\star$ can be given an directed tree structure with root $\mathcal{Q}_m$. Note that the Turing computation $\{\mathcal{Q}_0, \ldots, \mathcal{Q}_m\}$ is a path in this tree from the initial to the halting configuration.

**Definition 4.5** For each $\mathcal{Q} \in Cfg_\star$ we define three elements $\beta^j = \beta^j(\mathcal{Q})$ ($j \in \{0,1,2\}$) in $\mathbf{A}$ (they are also in $\mathbf{K}$). Assume that $\mathcal{Q} = \langle t, n, \mu_i \rangle$ and that $t(n) = r \in \{0,1\}$.

$$\beta^j(x) = \begin{cases} P_j & \text{if } x = s_0 \\ C_{ir}^{t(x)} & \text{if } x \in [b, n-1] \\ M_i^r & \text{if } x = n \\ D_{ir}^{t(x)} & \text{if } x \in [n+1, c] \end{cases}$$

Note that $q^j$ defined earlier is identical with $\beta^j(\mathcal{Q}_0)$.

Let $\Lambda = \{\beta^j(\mathcal{Q}) : \mathcal{Q} \in Cfg_\star \text{ and } j \in \{0,1,2\}\}$; also, let $K_0$ be the set of all $f \in K$ such that $0 \in f(N)$.

**Lemma 4.1** $K = K_0 \cup S \cup \Lambda$.

Proof: $\Lambda$ is generated from $S \cup \{q^0, q^1, q^2\}$ by applying just the machine operations. Thus $K_0 \cup S \cup \Lambda \subseteq K$.

For the reverse inclusion to fail, there must be some element $f \in K - K_0 - S - \Lambda$. In this case, $f = F(a_0, \ldots, a_{n-1})$ for some term $F(x_0, \ldots, x_{n-1})$ with each $a_i \in S \cup \{q^0, q^1, q^2\}$. Let this $F(\overline{x})$ be a term of minimal degree. Then we can write $F(\overline{x})$ as the composition of a basic operation $H$ with terms $F_i(\overline{x})$ of lesser degree. Then each $F_i(\overline{a}) \in K_0 \cup S \cup \Lambda$.

Case 1: $H \in \{\wedge, \wedge^i, T_0, T_1\} \cup \mathcal{M} \cup \mathcal{U}$. These operations all have 0 as their absorbing element in $\mathbf{B}(\mathcal{T})$, so each $a_i$ is in $S$ or $\Lambda$. We can immediately rule out $\wedge$, as $a_0 \wedge a_1 \notin K_0$ clearly implies that it's in $S \cup \Lambda$. We can similarly rule out $\wedge^0, \wedge^1, \wedge^2, T_0$ and $T_1$. We can also rule out the operations in $\mathcal{M}$ because of how $\Lambda$ was constructed. The operations in $\mathcal{U}$ must reduce to their corresponding operations in $\mathcal{M}$, as otherwise would require the existence of extra sequential type elements not in $S$. But that is impossible.

Case 2: $H \in \{J', S_1, S_2\}$. Note $J'(a_0, a_1, a_2) \leq a_0$, so $a_0$ must be in $S$ or $\Lambda$. But then it is clear to see that $f = J'(a_0, a_1, a_2)$ is also in $S$ or $\Lambda$. If $S_1(a_0, a_1, a_2, a_3) \notin K_0$ then $a_1, a_2$, and $a_3$ must be totally nonzero. If $a_1 \in S$ then we have that $a_1 = a_2 = a_3$ and $f = S_1(a_0, a_1, a_2, a_3) = a_1 \in S$. If $a_1 \in \Lambda$ then it must be the case that each $a_i \in \{\beta^j(\mathcal{Q}) : j \in \{0,1,2\}\}$ for some $\mathcal{Q} \in Cfg_\star$. But then $f$ is also an encoded configuration from that same set. The same argument works for $S_2$ and the result follows. ●

**Lemma 4.2** 1. Let $\mathcal{Q}, \mathcal{Q}' \in Cfg_\star$. Then there exists a polynomial function $F(x)$ of $\mathbf{K}$ satisfying $F(\beta^j(\mathcal{Q})) = \beta^j(\mathcal{Q}')$ for all $j \in \{0,1,2\}$.

2. Let $f \in S \cup \Lambda$. There is a polynomial function $F(x)$ of **K** such that $F(f) = q^2$ and $F(K_0) \subseteq K_0$.

3. If $F(x)$ is any nonconstant polynomial operation of **K** then $F(K_0) \subseteq K_0$.

4. Any polynomial of **K** that restricts to a subset of $\{q^0, q^1, q^2\}$ is order-preserving.

5. There is a polynomial function of **K** that when restricted to $\{q^0, q^2\}$ is the semi-lattice join operation.

Proof: The polynomials for (1) and (2) are built from machine operations. (1) essentially states that the tree $Cfg_\star$ is connected. For (2), we can find a machine operation that maps a sequential element of $S$ to a configuration element of $\Lambda$. By (1), we can map any element of $\Lambda$ to $\beta^j(\mathcal{Q}_0)$ for some $j \in \{0, 1, 2\}$. That is, to $q^j$. Then taking the meet with $q^2$ gives us $q^2$. Thus we can compose these operations to get a polynomial $F(x)$ of **K**. Of course, $F(K_0) \subseteq K_0$, as 0 is an absorbing element for the basic operations composed to make $F(x)$.

For (3) we consider a minimal counterexample to the claim. Let $F(x)$ be a unary polynomial of **K** of minimal degree among those for which $F(K_0) \not\subseteq K_0$. Then we can write $F(x) = H(F_0(x), \ldots, F_{n-1}(x))$, where $H$ is a basic operation of **K** and the $F_i(x)$ are polynomials of **K** of lesser degree. Then, each $F_i(x)$ is either constant, or $F_i(K_0) \subseteq K_0$. At least one must be non-constant. We can immediately rule out the possibility that $H$ is any operation besides $J'$, $S_1$, or $S_2$, as for all the other operations, 0 acts as an absorbing element of $\mathbf{B}(\mathcal{T})$.

Note that $J'(x, y, z)$ is 0-absorbing in $x$ and $y$, so if $H = J'$, then $F_0(x)$ and $F_1(x)$ must be constant. Furthermore, $J'(F_0, F_1, F_2(x)) \leq F_0$. Therefore, if $F_0 \in S$, then we must have $F_1 = F_0$ so $F(x) = J'(F_0, F_0, F_2(x)) = F_0 \wedge F_2(x)$, a contradiction, as this reduces to the case where $H = \wedge$, which we've already ruled out. If $F_0 \in \Lambda$, then $F_1 \in \Lambda$ as well, and both $F_0$ and $F_1$ agree on $N$. But in that case, we would also have that $F(x) = J'(F_0, F_1, F_2(x)) = F_0 \wedge F_2(x)$, again a contradiction.

If $F(x) = S_1(F_0(x), F_1(x), F_2(x), F_3(x))$ then $F_1(x), F(2(x)$ and $F_3(x)$ must all be totally nonzero, hence constants all in $S$ or $\Lambda$. If they are all in $S$, then $F_0(x)$ must also be totally nonzero, therefore constant, so $F(x)$ is constant as well. So we must have $F_1, F_2, F_3 \in \Lambda$, and there must be an element $f$ of $K_0$ such that $f(N) \subseteq \{1, 2\}$. But this is impossible, as such an element can only be generated from $S$, and it's simple to see that this cannot happen. If $H = S_2$, then by a similar argument we must have $u, v \in K_0$ such that for each $t \in N$, $u(t) = \overline{v(t)}$. But this should also be impossible. Elements in $\overline{V}$ can only be produced by using operations in $\mathcal{U}$, and then we can get barred elements in coordinates in $N$ only if we have $x, y, z \in K$ such that $x|_N \prec y|_N$ and $x|_N \prec z|_N$ with $y|_N \neq z|_N$ or else $x|_N \prec z|_N$, $y|_N \prec z|_N$ and $x|_N \neq y|_N$. But again, this is impossible, as the elements of $S|_N$ are a sequential set. Thus we have (3).

For (4), it should now be obvious, considering Lemma 4.1, that in **K**, the result of applying either $S_1$ or $S_2$ always lies in $K_0$. So the argument boils down to noticing that in $\mathbf{B}(\mathcal{T})$, $\wedge$, $T_0(x)$ and $T_1(x,y)$ are the only operations that are nonzero, nonconstant, and not the identity on $\{P_0, P_1, P_2\}$. These three operations are order preserving on $\{P_0, P_1, P_2\}$.

Finally for (5), Let $G$ be the composition of the Turing operations that simulate the Turing computation from $\mathcal{Q}_0$ to the halting state $\mathcal{Q}_m$. Then $G^{-1}(T_0(G(T_1(x,y))))$ restricts to a semilattice join operation. The polynomial $G$ allows us to compose $T_1$, which restricts to those elements of **K** that encode the initial configuration, with $T_0$, which restricts to those elements of **K** that encode the halting state. ●

**Lemma 4.3** *Let $f, g \in K$.*

1. $(f, g) \in \theta$ iff $f = g$ or $\{f, g\} \subseteq K_0$.

2. $(f, g) \in \overline{\theta}$ iff $(f, g) \in \theta$ or $\{f, g\} \subseteq \{\beta^0(\mathcal{Q}), \beta^1(\mathcal{Q}), \beta^2(\mathcal{Q})\}$ for some $\mathcal{Q} \in Cfg_\star$.

3. *The congruence $\theta$ has a unique cover in $\mathbf{ConK}$, namely $\overline{\theta}$.*

Proof: By part (3) of the lemma 4.2, all elements of $K_0$ are $\theta$-related. By part (2), no element of $K_0$ is $\theta$-related to any element in $K - K_0$. Now suppose that $\{f, g\} \subseteq K - K_0$ and that $f \wedge g \in K_0$. Then by part (2) of the lemma 4.2, there's a polynomial $F(x)$ of $K$ such that $F(f) = q^2$ but $F(f \wedge g) \in K_0$. Therefore, the polynomial $F(f \wedge x)$ shows that $(f, g) \notin \theta$. The final case to consider for part (1) is if $\{f, g\} \subseteq K - K_0$ but $f \wedge g \notin K_0$. In this case, it must be that $\{f, g\} \subseteq \{\beta^0(\mathcal{Q}), \beta^1(\mathcal{Q}), \beta^2(\mathcal{Q})\}$ for some $\mathcal{Q} \in Cfg_\star$. Then by part (1) of the lemma 4.2, we have $\{F(f), F(g)\} \subseteq \{q^0, q^1, q^2\}$ and $F(f) \neq F(g)$. If either $F(f)$ or $F(g)$ is $q^2$, then we have shown that $f$ and $g$ are not $\theta$-related. Otherwise $\{F(f), F(g)\} = \{q^0, q^1\}$ and we can use the polynomial $P(x) = q^0 \wedge F(x)$ for this case.

Suppose that $(f, g) \in K^2 - \theta$, and let $\theta'$ be the congruence generated by $\theta$ and $(f, g)$. We now show that $\theta'$ contains every pair $(\beta^j(\mathcal{Q}), \beta^{j'}(\mathcal{Q}))$ with $\mathcal{Q} \in Cfg_\star$ and $\{j, j'\} \subseteq \{0, 1, 2\}$. From the definition of $\theta$, it follows that $(q^2, h) \in \theta'$ for some $h \neq q^2$. If $h \in K_0$, then it follows from the previous lemma that all elements of $\Lambda$ are collapsed to $0/\theta$ by $\theta'$. This is the case even if $q^2 \wedge h \in K_0$. If $h = q^0$, then we use $T_1(q^2, x)$ to map $(q^0, q^2)$ to $(q^1, q^2)$. If $h = q^1$, then with $G(x)$ as before, $G^{-1}(T_0(G((x)))$ maps $(q^1, q^2)$ to $(q^0, q^2)$. Thus, $\{q^0, q^1, q^2\}$ is contained in one $\theta'$-class, and by lemma 4.1(1), the result follows. (3) follows from this last paragraph, as the pair $\{(q^0, q^2)\} \cup \theta \subseteq \theta'$ for each such $\theta'$. Part (2) can be concluded by showing that the equivalence relation $\Gamma$ over $K$ defined by

$$f \Gamma g \iff f = g \text{ or } \{f, g\} \subseteq K_0 \text{ or } \{f, g\} = \{\beta^j(\mathcal{Q}), \beta^{j'}(\mathcal{Q})\} \text{ for some } j, j', \mathcal{Q}$$

is a congruence of $K$. This is easy to check by showing that $\Gamma$ is preserved by every basic operation of $K$. We take care to mention that the reason we require $\{-1, 0, 1\} \subseteq N$ is that otherwise, $\{q^0, q^1, q^2\}$ would lack either a coordinate with $C_{10}^0$ or $D_{10}^0$ or $M_1^0$, and then $\Gamma$ would not be closed under $\wedge^0, \wedge^1, \wedge^2$. These three special operations will be used in the next section. •

**Lemma 4.4** *The type of the prime congruence $(\theta, \overline{\theta})$ in $\mathbf{ConK}$ is 4.*

Proof: By lemma 4.3, the pair $(q^0, q^2)$ belongs to $\overline{\theta} - \theta$. The induced algebra $\mathbf{K}|_{\{q^0, q^2\}}$ includes the meet operation $\wedge|_{\{q^0, q^2\}}$, the join operation from lemma 4.2(5), and by Lemma 4.2(4), the polynomial clone on $\mathbf{K}|_{\{q^0, q^2\}}$ is order-preserving. Hence $\mathbf{K}|_{\{q^0, q^2\}}$ is a minimal algebra of type 4. To show that $\{q^0, q^2\}$ is actually a trace, consider the polynomial $e(x) = q^0 \wedge x$. Then $e(K) = \{q^0, q^2\} \cup T$ where $T$ is some subset of $K_0$. As $e(K)$ consists of two $\overline{\theta}$ classes, one of which splits into distinct $\theta$ classes, $e(x)$ does not collapse $\overline{\theta}$ to $\theta$, and $e(K)$ contains a $(\overline{\theta}, \theta)$-trace. The only possibility is that this trace is $\{q^0, q^2\}$. •

## 5 If $\mathcal{T}$ doesn't Halt

**Theorem 5.1** *If $\mathcal{T}$ doesn't halt then $4 \notin typ\{V(\mathbf{B}(\mathcal{T}))\}$.*

We assume contrary to this claim, that $4 \in typ\{V(\mathbf{B}(\mathcal{T}))\}$ but that $\mathcal{T}$ does not halt. In this situation, $V(\mathbf{B}(\mathcal{T}))$ will contain a finite subdirectly irreducible algebra $\mathbf{K}$ with type 4 monolith. Since $\mathbf{K}$ is finite in $HSP(\mathbf{B}(\mathcal{T}))$, we can assume that $\mathbf{K} \cong \mathbf{A}/\theta$, where $\mathbf{A} \subseteq \mathbf{B}(\mathcal{T})^{\mathbf{Y}}$, and $Y$ is a set of least finite cardinality for this to be possible. Since we are assuming that $\mathbf{K}$ is a subdirectly irreducible with type 4 monolith, then $\theta \in Con\mathbf{A}$ has a unique cover $\overline{\theta}$, and the type of $(\theta, \overline{\theta})$ is 4.

In the next lemma, we collect some basic facts about (type 4) traces that could occur in the algebra $\mathbf{A}$. The proof of these useful facts only use that $\mathbf{A}$ has a semi-lattice reduct, and that the indexing set $Y$ is of minimal cardinality.

**Lemma 5.1** *There exists elements $\beta, \beta' \in A$ such that $\beta' < \beta$ and for all $f, g \in A$:*

1. *The pair $(\beta', \beta)$ belongs to $\overline{\theta} - \theta$.*

2. *$f\theta\beta$ implies $f \geq \beta$.*

3. *$f < \beta, g < \beta$ and $(f, g) \in \overline{\theta}$ imply $(f, g) \in \theta$.*

4. *The following are equivalent:*

    *(a) $(f, g) \in \theta$;*

(b) For every polynomial $F$ of $\mathbf{A}$,
$$F(f) \wedge \beta = \beta \Leftrightarrow F(g) \wedge \beta = \beta.$$

5. For every $t \in Y$ there exists $h \in A$ such that $h < \beta$ and $h(s) = \beta(s)$ for all $s \in Y - \{t\}$.

6. The set $\{\beta', \beta\}$ is a trace for the congruence quotient $(\theta, \overline{\theta})$.

Proof: We choose $\beta$ to be minimal in the order $\leq$ of the semilattice $\langle A, \wedge \rangle$, among all those $x \in A$ such that $(x, y) \in \overline{\theta} - \theta$ for some $y < x$; we choose $\beta' < \beta$ to satisfy $(\beta', \beta) \in \overline{\theta} - \theta$. Then (i),(ii) and (iii) are consequences of semilattice calculations. For (2), If $f\theta\beta$ then $f \wedge \beta\theta\beta \wedge \beta = \beta$. Also, $f \wedge \beta'\theta\beta \wedge \beta' = \beta'$. Therefore $((f \wedge \beta), (f \wedge \beta')) \in \overline{\theta} - \theta$. By minimality of $\beta$, we have $f \wedge \beta = \beta$. That is, $f \geq \beta$. For (3), if $f, g < \beta$ and $(f, g) \in \overline{\theta} - \theta$ then both $(f, f \wedge g)$ and $(f \wedge g, g)$ belong to $\overline{\theta}$, but not both of them are in $\theta$. The pair that does not belong to $\theta$ contradicts the minimality of $\beta$.

For (4): Let $\gamma$ be the relation on $\mathbf{A}$ defined by $(x, y) \in \gamma$ iff $x, y \in A$ and for every unary polynomial function $F$ of $\mathbf{A}$, we have $F(x) \wedge \beta = \beta$ iff $F(y) \wedge \beta = \beta$. Clearly $\gamma$ is a congruence of $\mathbf{A}$. To show that $\theta \subseteq \gamma$, suppose that $(f, g) \in \theta$ and $F(x)$ is a polynomial of $\mathbf{A}$. If $F(f) \wedge \beta = \beta$, then $(F(g) \wedge \beta)\theta(F(f) \wedge \beta) = \beta$. By (2) then, $F(g) \wedge \beta \geq \beta$, so we have $F(g) \wedge \beta = \beta$. To show that $\gamma \subseteq \theta$, note that otherwise, $\gamma$ contains $\overline{\theta}$. But this is not the case, as taking $(\beta', \beta) \in \overline{\theta}$ and $F$ the identity, $F(\beta) \wedge \beta = \beta \neq \beta' = F(\beta') \wedge \beta$.

As for (5), Let $t \in Y$ and $p'_t$ be the projection of $\mathbf{A}$ to $\mathbf{B}(\mathcal{T})^{Y-\{t\}}$. Since $|Y|$ is minimal, there must be $(x_0, x_1) \in A^2 - \theta$ such that $p'_t(x_0) = p'_t(x_1)$. Since $x_0$ and $x_1$ are not $\theta$ related, (4) guarantees a polynomial $F$ of $\mathbf{A}$ such that (possibly switching $x_0$ and $x_1$) $F(x_0) \wedge \beta = \beta$ while $F(x_1) \wedge \beta < \beta$. Let $h = F(x_1) \wedge \beta$. Then $h < \beta$ and for $s \neq t, h(s) = \beta(s)$.

Finally, for (6), consider the polynomial function $F(x) = x \wedge \beta$. Since $F(x)$ does not collapse $\overline{\theta}$ to $\theta$, the image $F(A) = [0, \beta]$ must contain a $(\theta, \overline{\theta})$-trace $\{x, y\}$. By (3), $\beta$ must belong to this trace. Thus we can write this trace as $\{h, \beta\}$ with $h < \beta$. Now change notation so $h$ becomes $\beta'$, and notice that (1)-(6) remain true. •

**Corollary**: $\beta$ is a totally non-zero function, and $\{\beta', \beta\}$ is a $\langle \theta, \overline{\theta} \rangle$-subtrace. Define $Y_0 = \beta^{-1}(\{P_0, P_1, P_2\})$ and $Y_1 = Y - Y_0 = \beta^{-1}(U \cup V \cup \overline{V})$.

**Lemma 5.2** *The following statements are true:*

1. $|Y| > 1$.

2. $|Y_0| \geq 1$.

3. $|Y_1| \geq 1$.

4. For $t \in Y_0$, $\beta'(t) \subseteq P$.

5. For $t \in Y_1$, $\beta'(t) = \beta(t)$.

Proof: (1)If $Y = 1$, then we can identify **A** with a subalgebra of $\mathbf{B}(\mathcal{T})$. We must analyze the polynomial structure of sets $\{\beta', \beta\} \subseteq \mathbf{B}(\mathcal{T})$ with $\beta' < \beta$.

Our analysis is simplified by the following claims: First, in $\mathbf{B}(\mathcal{T})$, if $\{a, b\}$ is a two element subset of $P \cup \{0\}$, then every polynomial image of $\{a, b\}$ is either also in $P \cup \{0\}$, or else is collapsed to a point. Thus, there is essentially no interaction between the polynomial structures of $\{0\} \cup P$ and $\{0\} \cup U \cup V \cup \overline{V}$.

Second, if $u \in \{1, 2\}$ and $v \in V \cup \overline{V}$, then the polynomials $f(x) = S_1(x, v, v, v)$ and $g(x) = S_2(x, \overline{v}, u, u, u)$ provide a polynomial isomorphism between $\{0, u\}$ and $\{0, v\}$. Therefore, we must be careful when calculating the polynomial structure of the two element ordered subsets of $\{0\} \cup U \cup V \cup \overline{V}$.

Let us start with $\{\beta', \beta\} \subseteq \{0\} \cup P$. In this case, it seems most efficient to perform the following analysis on the two element subsets of $\{0\} \cup P$ as in [5]:

Let's view these two element subsets as vertices of a directed graph, with an edge from $\{a, b\}$ to $\{c, d\}$ if and only if there exists $f \in Pol_1 \mathbf{B}(\mathcal{T})$ such that $\{c, d\} = \{f(a), f(b)\}$. It is easy to calculate these edges by using the fundatmental operations and then taking the transitive closure.

The strongly connected components of this graph are $M = \{\{P_0, P_1\}, \{P_2, P_0\}, \{P_2, P_1\}\}$ and $N = \{\{0, P_0\}, \{0, P_1\}, \{0, P_2\}\}$. If $C$ and $D$ are strongly connected components, then we write $C < D$ to indicate that there is an edge from each member of $D$ to each member of $C$, but none back. We have that $N < M$. Then it is immediately clear that the only $N$ contains subtraces. That's because for any pair $\{\mu, \nu\}$ in $M$, $\langle \mu, \nu \rangle$ is in the congruence generated by its shadow $N$, and thus is not a subtrace.

The three pairs in $N$ are subtraces, are all polynomially isomorphic and have a type 5 polynomial structure.

The pairs in $M$, while not subtraces, have a type 3 polynomial structure because of the ternary discriminator term $t(x, y, z)$ induced by $S_1$ and $S_2$. This will be used in the following arguments.

We proceed to calculate the polynomial structure of the two element ordered sets in $\{0\} \cup U \cup V \cup \overline{V}$. Suppose that $F(x) \in Pol_1 \mathbf{A}$ is nonconstant on $\{0, u\}$, where $u \in U \cup V \cup \overline{V}$. Then we claim that $F(0) = 0$. We give the typical argument: Consider a counterexample $F(x)$ of least degree. Then we can write $F(x) = H(F_1(x), \ldots, F_n(x))$, where $H$ is a basic operation of $\mathbf{B}(\mathcal{T})$, and the $F_i(x)$ are of lesser degree. Then for each $i$, $F_i(0) = 0$ or $F_i(x)$ is constant. $H$ must be one of $J'$, $S_1$, or $S_2$, as all the other basic operations are "0-absorbing". That is, their value is 0 whenever 0 is substituted for any variable. $H$ cannot be $S_1$ or $S_2$, as these operations are 0-absorbing too, when the range is restricted to $\{0\} \cup U \cup V \cup \overline{V}$. Now the last possibility, $J'(F_1(x), F_2(x), F_3(x))$, takes the value 0 when $F_1(x)$ or $F_2(x)$ is 0. Therefore, $F_1$ and $F_2$ must be constants, and in fact, $F_1 = \overline{F_2} \in V \cup \overline{V}$. But then $J'(F_1, \overline{F_1}, F_3(x)) = F_1$, a constant.

Therefore, none of these sets $\{0, u\}$ has a polynomial that switches the two elements. The semilattice meet $\wedge$ certainly restricts to $\{0, u\}$, so once we demon-strate that there is no pseudo-join, we will have shown that the polynomial structure is of type 4.

Let $u, v \in U \cup V \cup \overline{V}$. Then we say that the binary polynomial $F(x, y)$ is "join-like" for $u$ and $v$ if $F(0,0) = 0$, and $F(0, u) = F(u, 0) = F(u, u) = v$. Assume that $F(x, y)$ is join-like for some $u$ and $v$, and no other join-like binary polynomial has lesser degree.

Write $F(x, y) = H(F_1(x, y), \ldots, F_n(x, y))$ where $H$ is a basic operation, and the $F_i(x, y)$ are of lesser degree. Then none of the $F_i(x, y)$ are join-like, so restricting $x$ and $y$ to $\{0, u\}$, either $F_i(x, y)$ is constant, or $F_i(x, 0) = 0$ and $F_i(x, u)$ is constant, or $F_i(0, y) = 0$ and $F_i(u, y)$ is constant, or $F_i(0, 0) = F_i(u, 0) = F_i(0, u) = 0$.

If $H$ is any basic operation other than $S_1$, $S_2$ or $J'$, then $H$ is 0-absorbing. Therefore one of the $F_i$ must be join-like, contrary to assumption. If $H$ is either $S_1$ or $S_2$, then since $v \in U \cup V \cup \overline{V}$, $H$ actually turns out to be 0-absorbing ($S_1(u, x, y, z)$ is 0-absorbing in $x$, $y$, and $z$. $S_1(0, x, y, z) \neq 0$ only if $\{x, y, z\} \subseteq P$, giving $S_1$ a value in $P$. But above we showed that there is no non-constant map from $\{0, P_i\}$ to $\{0, v\}$ for any $v \in U \cup V \cup \overline{V}$. A similar argument works for $S_2$.) Finally, if $H = J'$, so that $F(x, y) = J'(F_1(x, y), F_2(x, y), F_3(x, y))$, then since $J'(x, y, z)$ is 0-absorbing in $x$ and $y$, we must have that $F_1(x, y)$ or $F_2(x, y)$ are either join-like (impossible by assumption) or constant. If they are constant, then either $F_1 = \overline{F_2}$ whence $F(x, y)$ is constant, or else $F_1 = F_2$, whence $F(x, y) = F_1 \wedge F_3(x, y)$. But then $F_3(x, y)$ must be join-like, a contradiction. Thus there are no join-like polynomials in $\mathbf{B}(\mathcal{T})$. In particular, there are no pseudojoins for $\{0, u\}$ for any $u \in U \cup V \cup \overline{V}$.

Now let's move on to item (2). Take a $t \in Y$ such that $\beta'(t) \neq \beta(t)$. Then, since $typ\{\beta', \beta\} = 4$, $\{\beta'(t), \beta(t)\}$ must have a polynomial structure in $\mathbf{B}(\mathcal{T})$ at least as rich. Refering back to the analysis of the two element subsets of $\mathbf{B}(\mathcal{T})$ above, we must have that $\{\beta'(t), \beta(t)\} \subseteq P$, and thus $|Y_0| \geq 1$.

(3) If $Y = Y_0$, then we have that $\beta'(Y), \beta(Y) \subseteq P$.

But in this case, both $S_1$ and $S_2$ induce a discriminator term that restricts to $\{\beta', \beta\}$, and $typ\{\beta', \beta\}$ would be 3. Thus $Y_1$ cannot be empty.

Finally, (4) and (5) are both easy and follow from the work in (1).

$\beta'(t) \in P$ for $t \in Y_0$, for otherwise, we would have that $\beta'(t) = 0$, and then $\{\beta'(t), \beta(t)\}$ does not have a rich enough polynomial structure to allow $\{\beta', \beta\}$ to be type 4.

For (5) if $\beta'(t) \neq \beta(t)$ for any $t \in Y_1$, then $\beta'(t) = 0$ while $\beta(t) \in U \cup V \cup \overline{V}$ and projecting down to that coordinate would induce a type 4 lattice structure where only a type 5 semilattice structure exists. Thus we must have that $\beta'(t) = \beta(t)$ for all $t \in Y_1$. ●

**Lemma 5.3** *The following items are some forbidden configurations.*

1. *In $\mathbf{A}$ there is no member $f$ satisfying $f(Y_1) \subseteq \{1, 2\}$;*

2. In **A** there are no two members $f, g$ satisfying $f(t) = \overline{g(t)} \in V \cup \overline{V}$ for all $t \in Y_1$;

3. If $f \in A$ and for all $t \in Y_1$ we have $\beta(t) = f(t)$ or $\beta(t) = \overline{f(t)} \in V \cup \overline{V}$ while $f(Y_0) \subseteq P$, then $f$ agrees with $\beta$ on $Y_1$.

Proof: For (1) and (2), such an occurence would induce a type 3 structure on $\{\beta', \beta\}$, using operations $S_1$ or $S_2$ respectively. For (3) Assume that such an $f \in A$ exists, and that there is a $t \in Y_1$ with $f(t) = \overline{\beta(t)}$. By Lemma 5.1(5), we can choose $h \in A$ with $h(t) < \beta(t)$ and $h(s) = \beta(s)$ for all $s \neq t$.

Let $F(x) = J'(\beta, f, x)$. Since $h = h \wedge \beta \overline{\theta} h \wedge \beta'$ and both $h \wedge \beta, h \wedge \beta' < \beta$, we have by lemma 5.1(3) that $h\theta h \wedge \beta'$. Therefore, $\beta = F(h)\theta F(h \wedge \beta')$ and by lemma 5.1(2), $F(h \wedge \beta') \geq \beta$. But in fact the opposite is true: $F(h \wedge \beta') = \beta'$. For if $s \in Y_0$, then $f(s) \subseteq P$ and $F(h \wedge \beta')(s) = J'(\beta, f, h \wedge \beta')(s) = \beta(s) \wedge h(s) \wedge \beta'(s) = \beta'(s)$;

If $s \in Y_1$ then either $\beta(s) = \overline{f(s)}$ so $F(h \wedge \beta')(s) = J'(\beta, f, \beta' \wedge h)(s) = \beta(s) = \beta'(s)$ or else $\beta(s) = f(s)$ so $F(h \wedge \beta')(s) = \beta(s) \wedge \beta'(s) \wedge h(s) = \beta'(s)$.

We have a contradiction, so the lemma is proved. ●

Definition: For $f \in A$ put $\Lambda(f) = \{g \in A : g|_{Y_1} = f|_{Y_1} \text{ and } g(Y_0) \subset P\}$.

Definition: $A_\star$ is the smallest set $X \subseteq A$ such that $\Lambda(\beta) \subseteq X$ and whenever $F$ is an operation in $\mathcal{M}$ and $\{f_0, f_1, f_2\} \subseteq A$ with $F(f_0, f_1, f_2) \in X$ then $\{f_0, f_1, f_2\} \subseteq X$.

**Lemma 5.4** *Here are some facts about $A_\star$:*

1. *For all $f \in A_\star$ either $f \in U^Y$ and $f(Y_0) \subseteq \{1\}$ or else $f(Y_0) \subseteq P$ and $f(Y_1) \subseteq V \cup \overline{V} \cup U$.*

2. *If $f \in A_\star$ and $f(Y_0) \subseteq P$ then $\Lambda(f) \subset A_\star$.*

3. *Let $f \in A_\star$, $g \in A$ and $g(Y_0) \subseteq P$. If for all $t \in Y_1$ either $g(t) = \overline{f(t)} \in V \cup \overline{V}$ or $g(t) = f(t)$ then $g$ and $f$ agree on $Y_1$.*

4. *Let $F(f_0, f_1, g) = h \in A_\star$ where $F \in \mathcal{M}$. Then $f_0 \prec f_1$; moreover, if $f \in A$ then $f_0 \prec f$ implies $f_1 = f$, while $f \prec f_1$ and $f(Y_0) \subseteq \{1\}$ imply that $f_0 = f$.*

Proof: Statements (1) and (2) should be obvious. Statement (3) follows from lemma 5.3 and the definitions of $A_\star$ and of the operations in $\mathcal{M}$. Statement (4) is easy, and uses the operations in $\mathcal{U}$. ●

**Lemma 5.5** *For every nonconstant polynomial function $F(x)$ of **A**, $F^{-1}(A_\star) \subseteq A_\star$.*

Proof: We assume that $F(x)$ is a non-constant polynomial of minimal degree such that $F^{-1}(A_\star) \not\subseteq A_\star$. Write $F$ as $H(F_0(x), \ldots, F_n(x))$, where $H$ is a basic operation of $A$, and each $F_i(x)$ is of lesser degree than the degree of $F$.

From the lemma 5.4, $A_\star$ is an up-set in the semilattice order of $\mathbf{A}$ and 0 is not in the range of any $f \in A_\star$. From this, it is clear that $H$ cannot be $\wedge$ or any $\wedge^i$ or $T_0$ or $T_1$. It follows from Lemma 5.3 (1) and (2) that the operations $S_1$ and $S_2$ have their range disjoint from $A_\star$ and so $H \in \{S_1, S_2\}$ is impossible.

If $J'(j, k, l) = f \in A_\star$ for some $j, k, l \subseteq A$, then since $J'(x, y, z) \leq x$, we already have $j \in A_\star$.

Now $f \in A_\star$, so by lemma 5.4, it is either in $U^Y$ or in $P^{Y_0} \times (U \cup V \cup \overline{V})^{Y_1}$.

In the first case, we must have that for each $t \in Y$, $j(t) = k(t) = l(t) = f(t)$ whence $j, k$, and $l$ are already in $A_\star$, so it must be the case that case two holds.

In the second case, we must have that for each $t \in Y_1$, $j(t) = k(t)$ or $j(t) = \overline{k(t)} \in V \cup \overline{V}$, and for each $t \in Y_0$, $\{j(t), k(t), l(t)\} \subseteq P$. It follows that there is a polynomial $G(x)$ built only out of machine operations and constants, such that $G(j) \in \Lambda(\beta)$. But then $G(k)$ in an element of $\mathbf{A}$ such that for each $t \in Y_1$, $G(k)(t) = \beta(t)$ or $G(k)(t) = \overline{\beta(t)}$, and for $t \in Y_0$, $G(k)(t) \subseteq P$. Then by lemma 5.3(3), $G(k)$ agrees with $\beta$ on $Y_1$, and therefore, by applying $G^{-1}$, $j$ and $k$ must agree on $Y_1$. Therefore, $J'(j, k, l) = j \wedge l$. Thus $H$ cannot be $J'$ either.

By the definition of $A_\star$, $H$ cannot be a machine operation. This along with Lemma 5.4(4) shows that $H \in \mathcal{U}$ is also impossible. Thus there is no such $F$ and $H$. •

**Lemma 5.6** $A_\star \neq \Lambda(\beta)$.

Proof: Assume otherwise, that $A_\star = \Lambda(\beta)$. Then for every operation $F \in \mathcal{M} \cup \mathcal{U}$, the range of $F$ in $\mathbf{A}$ is disjoint from $\Lambda(\beta)$.

We are led to consider two cases, and aim for a contradiction.

Case 1: $\beta(Y_1) \not\subseteq V_{10}^0$.

Then the range of $T_1$ is disjoint from $\Lambda(\beta)$. Likewise for each $\wedge^i$. By lemma 5.3, the same is true for $S_1$ and $S_2$. Thus the only basic operations that can have an element of $A_\star$ in their range are $\wedge$, $T_0$, and $J'$.

If $x \wedge y \in \Lambda(\beta)$ then we must have $x, y \in \Lambda(\beta)$. Likewise, if $T_0(x) \in \Lambda(\beta)$ then $x \in \Lambda(\beta)$ as well. By the same argument as in Lemma 5.5, $J'(x, y, z) \in \Lambda(\beta)$ only if $x = y$ on $Y_1$, so $J'(x, y, z) = x \wedge z$.

Therefore, if $F(x) \in Pol_1\mathbf{A}$ maps $\{\beta', \beta\}$ into $\Lambda(\beta)$, then $F(x)$ must agree on the set $\{\beta', \beta\}$ with some polynomial of the algebra $\langle A, \wedge, T_0 \rangle$. But projecting down onto any coordinate $t \in Y$, it is easy to compute that any two element subset $x < y$ of $\langle B(\mathcal{T}), \wedge, T_0 \rangle$ has only a type 5 structure.

Case 2: $\beta(Y_1) \subseteq V_{10}^0$.

In this case, $\beta(Y_1)$ must equal $V_{10}^0$, or else we could use one of the operations $\wedge^i$ to derive a contradiction. For example, if $M_1^0 \notin \beta(Y_1)$ then we use $\wedge^1$: Let $F(x) = \beta \wedge^1 x$ and $G(x) = x \wedge \beta'$. Then $F(\beta) = \beta$ and $G(\beta') = \beta'$ while $F(\beta') = G(\beta) \notin A_\star$ as it contains a 0 in its range. Therefore by Lemma 5.5, $\{\beta', \beta\}$ is not polynomially isomorphic to either of its images in $F$ or $G$. But at least one of these images must be in $\overline{\theta} - \theta$ and therefore must be a $(\theta, \overline{\theta})$ subtrace that is polynomially isomorphic to $\{\beta', \beta\}$. This is a contradiction.

Therefore $\beta(Y_1) = V_{10}^0$ and the range of each $\wedge^i$ is disjoint from $\Lambda(\beta)$. The only operations that can have range in $\{\beta', \beta\}$ are then $\wedge$, $T_1$ and $J'$. Again, $J'(x, y, z)$ must reduce to $x \wedge z$ so we merely have to calculate the polynomial structure of $\langle P, \wedge, T_1 \rangle$. But this is quite easy, and has only type 5 semilattice structures on its two element subsets $x < y$.

Thus the claim is proved. •

Definition: $\Sigma = A_\star \cap U^Y$ and $\Omega = A_\star - \Sigma$

**Lemma 5.7** $\Sigma$ *is a sequentiable set* $\{f_n : n \in N\}$ *for some convex set $N$ of integers. There is a connected set* $\Omega_c \subseteq Cfg_N(\mathcal{T})$ *such that the set* $\{v|_{Y_1} : v \in \Omega\}$ *can be written as* $\{\beta'(\mathcal{Q}) : \mathcal{Q} \in \Omega_c\}$ *where*

$$A_\star|_{Y_1} = \{f_n|_{Y_1} : n \in N\} \cup \{\beta'(\mathcal{Q}) : \mathcal{Q} \in \Omega_c\}$$

*is an algebraically coded computation of $\mathcal{T}$ in the sense defined [1, section 4].*

Proof: This argument gets slightly involved. Basically, we are reworking the arguments in Lemmas 5.8, 5.9, and 5.10 of [3]. We will spread the proof of this across the next few Lemmas, finishing after Lemma 5.10.

Definition: If $f \in \Sigma$ and $v \in \Omega$ then we write $f \sim v$ to indicate that for all $t \in Y_1$, we have:

$$f(t) = 1 \iff v(t) \in \{C_{ir}^s, \overline{C}_{ir}^s : i \leq k, \{r, s\} \subseteq \{0, 1\}\}$$

$$f(t) = H \iff v(t) \in \{M_i^r, \overline{M}_i^r : i \leq k, r \in \{0, 1\}\}$$

$$f(t) = 2 \iff v(t) \in \{D_{ir}^s, \overline{D}_{ir}^s : i \leq k, \{r, s\} \subseteq \{0, 1\}\}$$

If $v, v' \in \Omega$ then we write $v \sim v'$ if and only if there is an $f \in \Sigma$ such that $f \sim v$ and $f \sim v'$. If $x \sim y$ holds then we say that $x$ and $y$ **have the same pattern**.

For $f, g \in \Omega$, we write $f \smile g$ to denote that for each $t \in Y_1$, either $\{f(t), g(t)\} \subseteq V$ or $\{f(t), g(t)\} \subseteq \overline{V}$.

We define $\Omega_0 = \Lambda(\beta)$ and inductively we define $\Omega_{n+1}$ to be the set of all $u \in A_\star$ such that $u \in \Omega_n$ or $F(x, y, u) \in \Omega_n$ for some $F \in \mathcal{M}$ and $\{x, y\} \in A_\star$. Finally, for $n \geq 0$ we define $\Sigma_n$ as the set of all $f \in \Sigma$ such that $f \sim u$ for some $u \in \Omega_n$.

**Lemma 5.8** *We can say the following about the previous definitions:*

1. *If $\{f, f', u\} \subseteq A_\star$, $F \in \mathcal{L}$ and $F(f, f', u) = v \in A_\star$, then $\{f, f'\} \subseteq \Sigma$, $\{u, v\} \subseteq \Omega$, $f' \sim u$, $f \sim v$ and $u \smile v$.*

2. *If $\{f, f', v\} \subseteq A_\star$, $F \in \mathcal{L}^o$ and $F(f, f', v) = u \in A_\star$, then $\{f, f'\} \subseteq \Sigma$, $\{u, v\} \subseteq \Omega$, $f' \sim u$, $f \sim v$ and $u \smile v$.*

3. If $\{f, f', u\} \subseteq A_\star$, $F \in \mathcal{R}$ and $F(f, f', u) = v \in A_\star$, then $\{f, f'\} \subseteq \Sigma$, $\{u, v\} \subseteq \Omega$, $f \sim u$, $f' \sim v$ and $u \breve{\prec} v$.

4. If $\{f, f', v\} \subseteq A_\star$, $F \in \mathcal{R}^o$ and $F(f, f', v) = u \in A_\star$, then $\{f, f'\} \subseteq \Sigma$, $\{u, v\} \subseteq \Omega$, $f \sim u$, $f' \sim v$ and $u \breve{\prec} v$.

5. $\{\Omega_n\}$ is an increasing sequence of sets whose union is $\Omega$. $\{\Sigma_n\}$ is an increasing sequence of sets whose union is $\Sigma$.

6. $\Sigma$ is a nonvoid sequentiable set and every member of $\Omega$ has the same pattern as some member of $\Sigma$.

Proof: Statements (1)-(4) are immediate consequences of our definitions of the operations in $\mathcal{M}$ and of the relations $\sim$ and $\breve{\prec}$. Statement (5) is a corollary of our definitions of $A_\star$, $\Omega$ and $\Sigma$. Since $A_\star \neq \Lambda(\beta)$, it follows that $\Sigma \neq \emptyset$.

We show by induction on $n$ that $\Sigma$ is sequentiable. Clearly, $\Sigma_0 = \{f_0\}$ where $f_0 \sim \beta$ is sequentiable. Now suppose that $\Sigma_n$ is sequentiable, say

$$\Sigma_n = \{g_0, \ldots, g_b\}$$

where $g_0 \prec g_1 \prec \ldots \prec g_b$. Let $f \in \Sigma_{n+1} - \Sigma_n$. Then $f \sim u$ for some $u \in \Omega_{n+1} - \Omega_n$. We shall show that either $f \prec g_0$ or else $g_b \prec f$.

Suppose then, that $F(x, y, u) = v \in \Omega_n$ with $F \in \mathcal{L} \cup \mathcal{R}^o$ and $\{x, y\} \subseteq \Sigma$. Then by (1) or (3) we must have that $y \sim u$ as well. Therefore $f = y$. We also have that $x \sim v$, so $x \in \Sigma_n$. Suppose that $x = g_a$ for $a < b$. Then we see that $g_a \prec f$. It follows from Lemma 5.4(4) that $f = g_{a+1} \in \Sigma_n$, but we are sssuming that $f \in \Sigma_{n+1} - \Sigma_n$. Therefore $x$ must be $g_b$, and $g_b \prec f$ as claimed.

On the other hand, it's possible that $F \in \mathcal{L}^o \cup \mathcal{R}$. Then a similar argument gives $f \prec g_0$.

What if we have $f, f' \in \Sigma_{n+1} - \Sigma_n$ with $g_b \prec f$ and $g_b \prec f'$ (or else $f \prec g_0$ and $f' \prec g_0$)? Then we again appeal to lemma 5.4(4). There must be $v \in \Omega_n$, $u \in \Omega_{n+1} - \Omega_n$, $x, y \in \Sigma$ and $F \in \mathcal{L} \cup \mathcal{R}^o$ ($\in \mathcal{R} \cup \mathcal{L}^o$) such that $F(x, y, u) = v$ ($F(x, y, u) = v$) with $f \sim u$ ($f \sim u$). Here $x \sim v$ ($y \sim v$), so $x \in \Sigma_n$ ($y \in \Sigma_n$). Therefore $x$ ($y$) must be $g_b$ ($g_0$). Then by this same lemma, $f = f'$ Finally, it follows from lemma 5.3(1) that $f^{-1}(H) \neq \emptyset$ for $f \in \Sigma$. Therefore, $\Sigma_{n+1}$ is sequentiable.

Now, since $\Sigma$ is finite, there must be an $n$ for which $\Sigma_n = \Sigma$. So $\Sigma$ is also sequentiable. ●

As a result of this lemma, we can write $\Sigma$ uniquely as $\{f_0, f_1, \ldots, f_n\}$ for some $n$, with $f_i \prec f_{i+1}$ for $i < n - 1$. For each $i < n$ we define $X_i = f_i^{-1}(\{H\})$, and

$$X_L = \bigcap_{i<n} f_i^{-1}(\{1\})$$

$$X_R = \bigcap_{i<n} f_i^{-1}(\{2\})$$

Then $Y$ is the disjoint union of this partition.

Definition: $\Phi : V \cup \overline{V} \to V$ is the function that takes as its value the unbarred version of its argument. For example, $\Phi(M_1^0) = \Phi(\overline{M}_1^0) = M_1^0$.

**Lemma 5.9** *Let $q_1, q_2 \in \Omega$. There exists a polynomial $F(x)$ of $\mathbf{A}$ such that $F(\Lambda(q_1)) = \Lambda(q_2)$ and for all $f \in \Lambda(q_1)$, $F(f)|_{Y_0} = f|_{Y_0}$.*

Proof: Clearly, for $i \in \{1, 2\}$ we can build polynomials $F_i(x)$ using only operations in $\mathcal{M}$ and constants in $\Sigma$ such that $F_i(\Lambda(q_i)) = \Lambda(\beta)$ and for all $f \in \Lambda(q_i)$, $F_i(f)|_{Y_0} = f|_{Y_0}$.

Now we set $F(x) = F_2^{-1}(F_1(x))$. This function satisfies the lemma, and it is a polynomial, as $F_2^{-1}(x)$ can be built out of the machine operations that reverse the ones used to build $F_2(x)$. •

**Lemma 5.10** *For all $v \in \Omega$ we have $v \smile \beta$. Also, for each $0 \leq j < n$ and $v \in \Omega$, the function $\Phi(v(s))$ is constant as $s$ varies over the set $X_j$.*

Proof: We can prove these statements by induction.

For $v \in \Omega_0 = \Lambda(\beta)$, we certainly have that $v \smile \beta$.

Now suppose that this assertion holds for all $u \in \Omega_n$, and take $v \in \Omega_{n+1} - \Omega_n$. Then there is an $F \in \mathcal{M}$ and $x, y \in \Sigma$ and $u \in \Omega_n$ such that $F(x, y, v) = u$. Then $u \smile v$ follows from the definition of $F$. By Lemma 5.8(5), this argument suffices.

For the second statement, we first make the weaker claim that if $f_j \sim v$ for $v \in \Omega$, then $\Phi(v(s))$ is constant as $s$ varies over $X_j$.

Fix $v \in \Omega$ and take $f_j \in \Sigma$ such that $f_j \sim v$. We make some assumptions to make the argument clearer, but with slight modifications, a similar argument works in each case. We assume that $j \neq b$ (so that $f_{j+1} \in \Sigma$), and that there is a $u \in \Omega$ such that $L_{irt}(f_j, f_{j+1}, u) = v$, for some $0 < i \leq k$, and $r, t \in \{0, 1\}$. By the definition of $L_{irt}$, it is immediate that $\Phi(v(s)) = M_m^t$ for each $s \in X_j$. This gives us the weaker claim.

For the full claim, fix $j$ and take any $v \in \Omega$ for which $\Phi(v(s))$ is constant as $s$ varies over $X_j$. Notice that for any $u \in \Omega$, if there are $F \in \mathcal{M}$ and $x, y \in \Sigma$ so that either $F(x, y, u) = v$ or $F(x, y, v) = u$, then $x \prec y$ and $\{(x(t), y(t)) : t \in X_j\}$ is one of the following sets: $\{(1,1)\}$, $\{(H,1)\}$, $\{(2,H)\}$ or $\{(2,2)\}$. Therefore, by the definition of $F$, $\Phi(v(s))$ must also be constant as $s$ varies over $X_j$. The full claim then follows from the weaker claim, and the construction in Lemma 5.9. •

Now we can finally finish the proof of Lemma 5.7. Lemma 5.8(6) shows that $\Sigma$ is a nonvoid sequentiable set, that we have written as $\Sigma = \{f_0, \ldots, f_n\}$ with $f_i \prec f_{i+1}$ for each $0 \leq i < n$. Now define $N = \{0, \ldots, n\}$.

A configuration for $\mathcal{T}$ is a triple $\mathcal{Q} = \langle t, m, \gamma \rangle$ in which $t$ is a tape (that is, a function from $Z$ to $\{0, 1\}$), $n \in Z$ and $\gamma \in \{\mu_0, \ldots \mu_k\}$ is an internal state. For each element $v \in \Omega$, we define $\mathcal{Q}(v)$, the configuration for $\mathcal{T}$ encoded by $v$ as follows:

$$t(i) = \begin{cases} 0 & \text{for } i < 0 \\ k & \text{if } 0 \leq i \leq n \text{ and } \Phi(v(s)) = M_j^k \text{ or} \\ & \text{if } \Phi(v(s)) = C_{jr}^k \text{ or} \\ & \text{if } \Phi(v(s)) = D_{jr}^k \text{ as } s \text{ varies over } X_i. \\ 0 & \text{for } i > n \end{cases}$$

$$m = \text{ the unique } i \in N \text{ such that } f_i \sim v$$

$\gamma = $ the unique machine state $\mu_l$ such that $\Phi(v(s)) = M_l^k$ as $s$ varies over $X_m$.

We define $\Omega_c$ to be the collection of configurations encoded by the elements of $\Omega$ as defined above.

It is straightforward to check that if $u, u' \in \Omega$ encode the configurations $\mathcal{Q}, \mathcal{Q}'$, then $\mathcal{T}(\mathcal{Q}) = \mathcal{Q}'$, if and only if there is an $F \in \mathcal{L} \cup \mathcal{R}$ such that for some $f_i, f_{i+1} \in \Sigma$ we have that $F(f_i, f_{i+1}, u) \in \Lambda(u')$. Thus, it quickly follows from Lemma 5.9, that $\Omega_c$ is a connected set of configurations for $\mathcal{T}$.

Thus, $\Sigma|_{Y_1} \cup \Omega|_{Y_1}$ is really an algebraically encoded Turing computation of $\mathcal{T}$ as defined in [1, section 4]. $\bullet$

**Lemma 5.11** *Given $u \in \Omega$, there is at most one forward machine operation $F \in \mathcal{L} \cup \mathcal{R}$, such that $F(f_i, f_{i+1}, u) \in \Omega$.*

Proof: Suppose that $\mathcal{Q}(u) = \langle t, m, \mu_i \rangle$. Then there is at most one (none, if $\mu_i$ is the halting state $\mu_0$) command from $\mathcal{T}$ that begins with $\mu_i t(m)$. If the full command is $\mu_i t(m) s T \gamma$ where $s \in \{0,1\}$, $T \in \{L, R\}$ and $\gamma \in \{\mu_0, \ldots, \mu_k\}$ and if $t(m-1) = d$ if $T = L$ (or if $t(m+1) = d$ if $T = R$), then considering the definitions, $L_{it(m)d}$ (or $R_{it(m)d}$) is the unique forward command that gives a member of $\Omega$ in its output with $u$ as one of the arguments. $\bullet$

Corollary: The set of configurations encoded by $\Omega$ can be given a directed graph structure. We take the configurations as the vertices. There is an edge from $\mathcal{Q}$ to $\mathcal{Q}'$ exactly when $\mathcal{T}(\mathcal{Q}) = \mathcal{Q}'$. Likewise, using $\{\Lambda(v) : v \in \Omega\}$ as the vertices; with an edge from $\Lambda(u)$ to $\Lambda(v)$ if there is an operation $F \in \mathcal{R} \cup \mathcal{L}$ and $f_i, f_{i+1} \in \Sigma$ such that $F(f_i, f_{i+1}, u) \in \Lambda(v)$. By Lemma 5.9 this directed structure is connected. By this lemma, there is at most one edge leaving any vertex. Therefore, if there is a vertex with no outgoing edges, then we must have a directed tree structure terminating in this vertex.

**Lemma 5.12** *There cannot exist $q_1, q_0 \in \Omega$ such that $q_1(Y_1) \subseteq V_{10}^0$ and $q_0(Y_1) \subseteq V_0$.*

Proof: Suppose that such elements of $\Omega$ do exist. Then $q_1$ encodes an initial-state configuration $\mathcal{Q}_1$ and $q_0$ encodes a halting-state configuration $\mathcal{Q}_0$. By the corollary to the previous Lemma, we have a directed tree structure with $\mathcal{Q}_0$ as the terminal vertex. Therefore, there is a path from the initial configuration $\mathcal{Q}_1$ to the halting configuration. But this path is just a halting computation of $\mathcal{T}$ from a blank tape in machine state $\mu_1$, which contradicts our assumption about $\mathcal{T}$.  •

**Lemma 5.13** *There exists no $q \in \Omega$ with $q(Y_1) \subseteq V_0$.*

Proof: Otherwise, we could assume, without loss of generality, that $\beta'(Y_1), \beta(Y_1) \subseteq V_0$.

It follows from Lemma 5.5 that $\{\beta', \beta\}$ can only be polynomially equivalent to other pairs in $A_\star$, in fact in $\Omega$. Therefore, it is best to restate the claim of [1, Lemma 6.12]:

Claim: Let $F(x)$ be a polynomial function of $\mathbf{A}$ such that $\{F(\beta'), F(\beta)\}$ is a subset of $\Omega$. There exist polynomial functions $P(x)$, $Q(x)$ such that $F(x) = Q(P(x))$ for $x \in \Lambda(\beta)$ and $P|_{\Lambda(\beta)}$ is identical with a polynomial of the algebra $\langle \Lambda(\beta), \wedge, T_0 \rangle$ while $Q : \Lambda(\beta) \to \Lambda(F(\beta))$ is the unique polynomial function discussed in Lemma 5.9.

Proof of claim: see [1].

Therefore, by this claim, the problem of calculating the polynomial structure on $\{\beta', \beta\}$ reduces to calculating the polynomial structure on $\Lambda(\beta)$ generated by $\wedge$ and $T_0$. But this structure does not admit a lattice structure on a two element subset.  •

Now we can prove the theorem. We must modify the claim this time. As a result of Lemma 5.13, the operation $T_0$ has range disjoint from $\Omega$. Likewise for $S_1$, $S_2$, $\wedge^0$, $\wedge^1$, and $\wedge^2$. Thus, in the case that there is $q \in \Omega$ with $q(Y_1) \subseteq V_{10}^0$, then $P(x)$, when restricted to any coordinate in $Y_0$ must be identical with a polynomial in $\langle \{0, P_2, P_0, P_1\}, \wedge, T_1 \rangle$. But this polynomial structure is not rich enough to admit type 4. If there is $q \in \Omega$ such that $q(Y_1) \not\subseteq V_{10}^0$, then the polynomial structure is even smaller, hence no type 4 structure appears either.

Thus in the case that $\mathcal{T}$ does not halt, 4 does not appear in the typeset of $V(\mathbf{B}(\mathcal{T}))$.  •

[5] J. Berman and S. Seif *An approach to tame congruence theory via subtraces,* Algebra Universalis (1993), 479-520.